\documentclass[review,sort&compress]{elsarticle}

\usepackage{lineno,hyperref}
\modulolinenumbers[5]

\usepackage{tikz} 
\usetikzlibrary{shapes} 
\usepackage{pgfplots} 
\usepackage{amssymb,amsmath,amsthm} 
\usepackage[per-mode=symbol,exponent-product = \cdot]{siunitx} 
\usepackage[]{todonotes} 
\usepackage{mathtools}
\usepackage{algorithm,algpseudocode}
\usepackage{listings}

\newtheorem{remark}{Remark} 

\journal{Applied Mathematics and Computation}

\newcommand{\abs}[1]{\left|#1\right|}

\newcommand{\tab}{\hspace{0.5cm}}
\newcommand{\pf}[2]{\frac{\partial #1}{\partial #2}}

\newcommand{\ds}{\displaystyle}
\newcommand{\lr}[1]{\left( #1 \right)}
\newcommand{\LR}[1]{\left[ #1 \right]}
\newcommand{\Lr}[1]{\left\{ #1 \right\}}

\DeclareMathOperator{\sign}{sign}

\newcount\colveccount
\newcommand*\colvec[1]{
        \global\colveccount#1
        \begin{bmatrix}
        \colvecnext
}
\def\colvecnext#1{
        #1
        \global\advance\colveccount-1
        \ifnum\colveccount>0
                \\
                \expandafter\colvecnext
        \else
                \end{bmatrix}
        \fi
}

\newcommand{\Vb}{\mathbf{b}}
\newcommand{\Ve}{\mathbf{e}}
\newcommand{\Vr}{\mathbf{r}}

\newcommand{\Vm}{\mathbf{m}}
\newcommand{\Vn}{\mathbf{n}}

\graphicspath{{tikz-imgs/}} 

\begin{document}

\begin{frontmatter}

\title{Adaptive Refinement Strategies for the Simulation of Gas Flow in Networks using a Model Hierarchy}


\author[TUD]{Pia Domschke}
\ead{domschke@mathematik.tu-darmstadt.de}

\author[TUB]{Aseem Dua}
\ead{dua.aseem@gmail.com}

\author[TUB]{Jeroen J.\ Stolwijk\corref{corauth}}
\ead{stolwijk@math.tu-berlin.de}
\cortext[corauth]{Corresponding author}

\author[TUD,ESE]{Jens Lang}
\ead{lang@mathematik.tu-darmstadt.de}

\author[TUB]{Volker Mehrmann}
\ead{mehrmann@math.tu-berlin.de}

\address[TUD]{Department of Mathematics, TU Darmstadt, Dolivostr.\ 15, 64293~Darmstadt, Germany}
\address[ESE]{Graduate School of Energy Science and Engineering, TU Darmstadt, \\Jovanka-Bontschits-Str.\ 2, 64287 Darmstadt, Germany}
\address[TUB]{Institut f\"ur Mathematik, MA 4-5, TU Berlin, \\Stra{\ss}e des 17.\ Juni 136, 10623 Berlin, Germany}

\begin{abstract}
A model hierarchy that is based on the one-dimensional isothermal Euler equations of fluid dynamics is used for the simulation and optimisation of gas flow through a pipeline network. Adaptive refinement strategies have the aim of bringing the simulation error below a prescribed tolerance while
keeping the computational costs low. While spatial and temporal stepsize adaptivity is well studied in the literature, model adaptivity is a new field of research. The problem of finding an optimal refinement strategy that combines these three types of adaptivity is a generalisation of the unbounded knapsack problem.
A~refinement strategy that is currently used in gas flow simulation software is compared to two novel greedy-like strategies. Both a theoretical experiment and a realistic gas flow simulation show that the novel strategies significantly outperform the current refinement strategy with respect to the computational cost incurred.
\end{abstract}


\end{frontmatter}

\linenumbers

\section{Introduction}

The simulation of gas flows in pipeline networks is a topic of research that has
been studied at various scales: from the individual pipeline to the entire network. 
Studies in control and optimisation of gas supply in a dynamic supply-demand environment
strongly depend on large scale simulations of pipeline networks. In the last decades, considerable research on the modelling, simulation and optimisation of gas flow through pipeline networks has been conducted, see e.g.~\cite{AdeZ1995,BanHK2006-2,BanHK2006,ChaKW2005,ColG2006,DomKL2015,EhrS2005,EhrS2004,HerMS2010,
KeT2000,MarMM2006,Osi2001,Ste2007}.
Depending upon requirement, there exist
multiple models to predict the system behaviour with varying levels of accuracy. Generally, more accurate models  are computationally more expensive. Hence, in order to make real-time decisions, an appropriate trade-off between accuracy and computational complexity should be made.
This can be achieved by using a hierarchy of models, where the models can be adaptively switched during the simulation process. Beside the models, the discretisation mesh may be varied in space and time, which places
the demand for a strategy to automatically steer the simulation by changing
the models and the discretisation meshes. This steering is based on
simulation error estimates, which have been studied in detail in~\cite{DomschkePhD}.

Since the simulations are the basis for decisions in the optimisation and control of the gas flow, the reliability of the simulation is of prime importance. The simulation is to be carried out such that the
relative error in the state or in a functional of interest is below a specified tolerance.
Starting with a coarse simulation, an adaptive strategy is used to bring the error below the tolerance by refining the discretisation in time and space or refining models, i.e., shifting to a model of higher accuracy. Hence, we have three different refinement possibilities for each pipe $j \in \mathcal{J}_p$ of the pipeline network, where $\mathcal{J}_p$ denotes the set of pipes in the network.
These refinement possibilities are indexed by $i = 1,\ldots,3 N_p$, where $N_p \coloneqq \abs{\mathcal{J}_p}$ is the number of pipes. Refinements are to be chosen
such that the computational costs are kept low. 
We define an optimal refinement strategy as a strategy which returns the solution of the constrained optimisation problem
\begin{equation}
\begin{aligned}\label{eq:gen_Knapsack}
\min_{r_i} \tab &c + \sum_{i=1}^{3 N_p} \sum_{k=1}^{r_i} v_{ik} \\
\text{s.t.} \tab &\eta - \sum_{i=1}^{3 N_p} \sum_{k=1}^{r_i}w_{ik} \leq \text{tol}_\eta.
\end{aligned}
\end{equation}
Here, for each refinement 
possibility $i$, $r_{i}$ is the number of refinements, $w_{ik}$ is the relative error reduction due to the $k^\text{th}$ subsequent refinement and $v_{ik}$ is the corresponding cost addition. The constants $c$ and $\eta$ denote the cost and relative error of the starting simulation, respectively. We note that if $v_{ik}$ and~$w_{ik}$ are constant for all $k$, then this problem is equivalent to the unbounded knapsack problem which is NP-hard, see e.g.~\cite{KelPP04}. 
In this paper we aim to find a good approximation to the solution of this generalisation of the knapsack problem~\eqref{eq:gen_Knapsack}.
For this, we propose three
adaptive refinement strategies
which return approximate solutions of \eqref{eq:gen_Knapsack}.
The ideas presented here are for the example of pipeline networks. By generalising pipes to functional sub-domains, the principles of adaptive refinement can be extended to simulations for other applications which use a model hierarchy,
e.g., power grids and water supply networks.

This paper is organised as follows. Section \ref{sec:Model_Hierarchical_Simulations_of_Gas_Flows} introduces a model hierarchy for the simulation of gas flow through a single pipe as well as the aim of the adaptive refinement strategies. Section \ref{sec:strategies} describes the three proposed refinement strategies and Section \ref{sec:des_of_expt} introduces the design of a synthetic experiment. The results of this synthetic experiment are contained in Section \ref{sec:results} and an application of the refinement strategies to a realistic gas network simulation is given in Section~\ref{Application_to_a_Network_Simulation}. Finally, some conclusions are contained in Section \ref{Conclusions}.

\section{Simulation of Gas Flow}
\label{sec:Model_Hierarchical_Simulations_of_Gas_Flows}
In this section, we give an example of a model hierarchy for gas flow simulations. We then outline a framework for the adaptive simulation of pipeline networks using the given model hierarchy. Finally, we highlight the aims of the refinement strategies in the simulation framework.

\subsection{The Model Hierarchy}
\begin{figure}[!b]
\centering
\begin{tikzpicture}
\tikzstyle{model} = [rectangle, rounded corners = 5pt, minimum width=2cm, minimum height=1cm, text centered, draw=black,text width=4cm]
\tikzstyle{arrow} = [style=->,>=stealth]
\linespread{1.3}
\node (ISO1) [model] at (7,4.5) {Isothermal Euler Equations (M1)};
\node (ISO2) [model] at (7,2.25) {Semilinear Model (M2)};
\node (ISO4) [model] at (7,0) {Algebraic Model (M3)};
\draw[arrow] (ISO1.south) -| (ISO2) node [right, near end, fill=none] {$\pf{}{x}\lr{\rho v^2} = 0$};;
\draw[arrow] (ISO2.south) -| (ISO4) node [right, near end, fill=none] {$\pf{\rho}{t} = \pf{}{t}(\rho v) = h^\prime = 0$};;
\end{tikzpicture}
\caption{The model hierarchy that is considered in this paper.}
\label{ModelH}
\end{figure}
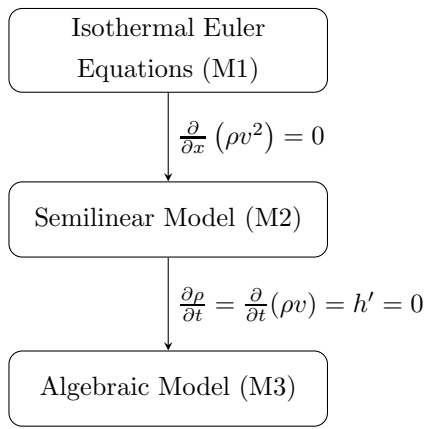

As an example, we take a three-model hierarchy for the gas flow simulations, discussed in detail in \citep{Domschke2011b}. At the highest level, we have the \emph{isothermal Euler equations} relating the gas density $\rho$ and gas flow $\rho v$ by
\begin{align}
\tag{\text{M}1}
\label{M1}
\begin{split}
\pf{\rho}{t} + \pf{}{x}\lr{\rho v} &= 0, \\
\pf{}{t}\lr{\rho v} + \pf{}{x}\lr{p+\rho v^2} &= -\frac{\lambda}{2D}\rho v\abs{v} - g\rho h',
\end{split}
\end{align}
together with the equation of state for real gases $p= \rho z(p)RT$ with compressibility factor $z(p) = 1 - \alpha p$ for a constant $\alpha \in \mathbb{R}^+$. Here, $p$ denotes the pressure, $v$ the velocity of the gas,
$R$ the specific gas constant, and $T$ the temperature. Further, $\lambda > 0 $ is the Darcy friction coefficient, $D$ the pipe diameter, $g = \SI{9.81}{\meter\per\second\squared}$ the acceleration due to gravity, and $h^\prime$ the slope of the pipeline. If the term $\pf{}{x}\lr{\rho v^2}$ is small, it can be dropped resulting in a \emph{semilinear model}
\begin{align}
\tag{\text{M}2}
\label{M2}
\begin{split}
\pf{\rho}{t} + \pf{}{x}\lr{\rho v} &= 0, \\
\pf{}{t}\lr{\rho v} + \pf{p}{x} &= -\frac{\lambda}{2D}\rho v\abs{v} - g\rho h'.
\end{split}
\end{align}
A further simplification of assuming a stationary state and zero slope $h^\prime = 0$ yields a system of two ODEs, which can be solved analytically and are referred to as the \emph{algebraic model}
\begin{align}
\tag{\text{M}3}
\label{M3}
\begin{split}
\rho v &= \text{constant}, \\
p(x) &= \sqrt{p_{\text{in}}^2 - \frac{\lambda c^2x}{D}\rho v \abs{\rho v}}.
\end{split}
\end{align}
Here, $c = \sqrt{p/\rho}$ denotes the speed of sound within the gas and $p_{\text{in}} = p(0)$ the input pressure.
The three models are shown in hierarchical form in Fig.\ \ref{ModelH}. The model hierarchy is set in the decreasing order of accuracy for our purpose of gas flows. Each pipe in the network is simulated using one of these three models and varying discretisation stepsizes in space and time.

\subsection{Adaptive Gas Network Simulation}
\begin{figure}[!bh]
\centering
\begin{tikzpicture}
\tikzstyle{startstop} = [rectangle, rounded corners, minimum width=3cm, minimum height=1cm,text centered, draw=black, scale = 0.9]
\tikzstyle{process} = [rectangle, minimum width=1cm, minimum height=1cm, text centered, draw=black,text width=4cm, scale = 0.9]
\tikzstyle{decision} = [diamond, minimum width=1cm, minimum height=1cm, text centered, draw=black,text width=2cm, scale = 0.9]
\tikzstyle{arrow} = [style=->,>=stealth]
\tikzset{	connector/.style={-latex,font=\scriptsize},
			rectangle connector/.style={connector,to path={(\tikztostart) -- ++(#1,0pt) \tikztonodes |- (\tikztotarget) }, pos=0.5},					 		 rectangle connector/.default=-2cm,
			 straight connector/.style={connector,to path=--(\tikztotarget) \tikztonodes }
		}
\linespread{1.3}
\node (start) [startstop] at (0,1.8) {start: $k=0$};
\node (simulate) [process,text width= 3cm] at (3.5,0){Simulate for $\LR{t_k, t_{k+1}} $};
\node (error) [process] at (8,0){A posteriori error estimates give $\Ve_m$, $\Ve_x$, $\Ve_t$};
\node (error check) [decision] at (8,-2.7){Total error $< \text{tol}$?};
\node (adaptation) [process] at (3.5, -2.7){Adaptive strategy returns refinements $\Vr_m$, $\Vr_x$, $\Vr_t$};
\node (conditions)[process,text width = 3cm ]at (0,0){$\Vm_k, \Vn_{x,k}, \Vn_{t,k}$ $\Vr_m\!=\!\Vr_x\!=\!\Vr_t\!=\!0$};
\node (update) [process, text width=2cm] at (0,-2.7) {$k=k+1$};
\node (coarsen) [process] at (3.5, -4.7){Coarsen model and discretisation where possible};

\draw[arrow] (start.south) -| (conditions);
\draw[arrow] (conditions.east) |- (simulate.west);
\draw[arrow] (simulate.east) |- (error.west);
\draw[arrow]  (error.south) -| (error check.north);
\draw [straight connector] (error check) to node[above,midway] {No} (adaptation);
\draw [rectangle connector=0cm] (error check.south) to node[above, midway, xshift=-1.5cm,yshift=-0.5cm]{Yes} (coarsen.east);
\draw[arrow]  (adaptation.north) -| (simulate.south);
\draw[arrow] (update.north) -| (conditions.south);
\draw[arrow] (coarsen.west) -| (update.south);
\end{tikzpicture}
\caption{Gas flow adaptive simulation process using an adaptive refinement strategy which returns a set of refinements $\Vr_m$, $\Vr_x$ and $\Vr_t$.}
\label{SimFlow}
\end{figure}
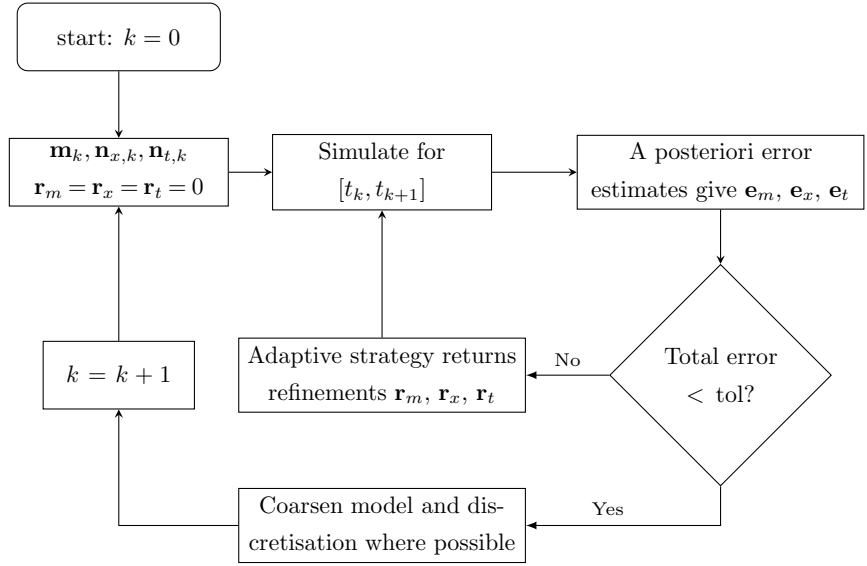

We consider a gas flow simulation over a pipeline network $\mathcal{J}_p$.
The simulation time $\LR{0,T}$ is divided into time intervals of equal size $\LR{t_{k}, t_{k+1}}$, \mbox{$k = 0, 1, \dots, N\!-\!1$} and $t_{N}  = T$. Given a starting model distribution over the network $\Vm_0 = \LR{m_{0,1}, m_{0,2}, \ldots, m_{0,N_p}}^T$, with $m_{0,j} \! \in \! \{1,2,3\}$, and a corresponding discretisation in space~$\Vn_{x,0}$ and in time~$\Vn_{t,0}$, a simulation is run for $\LR{t_0, t_{1}}$. We obtain error distributions along the network using \textit{a posteriori} error estimations (see \cite{DomschkePhD}): $\Ve_m = \LR{e_{m,1}, e_{m,2}, \ldots, e_{m,N_p}}^T$ for the model errors and $\Ve_x = \LR{e_{x,1}, e_{x,2}, \ldots, e_{x,N_p}}^T$ and $\Ve_t= \LR{e_{t,1}, e_{t,2}, \ldots, e_{t,N_p}}^T$ for the spatial and temporal discretisation errors, respectively.
These error estimators are derived using a dual weighted residual method with a user-defined functional of interest~$M(\mathbf{u})$, where $\mathbf{u} = \LR{\rho, \rho v}^T$, see \cite{DomschkePhD}.
The simulation error for a single pipe is the sum of all three errors. For the simulation to be valid, the relative error
must be below a given tolerance $\text{tol}$. Hence, we require that
\begin{equation}
\label{TolTarget}
\frac{\sum_{j \in \mathcal{J}_p} \lr{ e_{m,j} + e_{x,j} +e_{t,j} }}{\abs{M(\mathbf{u}^h)}} < \text{tol}.
\end{equation}
If the tolerance is not achieved, models 
and discretisation meshes are refined. The task of deciding the required refinements is made by an
adaptive strategy.
A switch to a higher model in the hierarchy is
called a \emph{model refinement} and a refinement of the mesh is called a \emph{discretisation refinement}. With the new models and discretisations we re-simulate for the time interval and continue the cycle. Once the solution meets the tolerance requirements, the models and discretisations are coarsened where possible, the simulation progresses to the next time interval and the cycle repeats. This simulation flow is shown in Figure~\ref{SimFlow} for an interval $\LR{t_k, t_{k+1}}$.

\subsection{Structure and Aim of Adaptive Strategies}

Our focus
lies on finding adaptive
strategies that control the errors and drive the simulation. For the time interval $\LR{t_k, t_{k+1}}$,
an adaptive
strategy takes as input the error distributions $\Ve_m, \Ve_x, \Ve_t$, the
model distribution~$\Vm_k$ and the number of nodes in the spatial and temporal discretisations $\Vn_{x,k}, \Vn_{t,k}$. The strategy
returns a refinement scheme $\Vr_m = \LR{r_{m,1},\ldots,r_{m,N_p}}^T$, $\Vr_x = \LR{r_{x,1},\ldots,r_{x,N_p}}^T$ and $\Vr_t = \LR{r_{t,1},\ldots,r_{t,N_p}}^T$, where $r_{m,j} \! \in \! \{0,1,2\}$ and $r_{x,j}, r_{t,j} \! \in \! \mathbb{N}$
denote the number of refinements to be made in the models and in the discretisations for all pipes $j \in \mathcal{J}_p$ such that constraint \eqref{TolTarget} is satisfied. The aim of the adaptive strategies is to achieve this constraint while
keeping the computational costs that are incurred in the simulation low. 
%
%

\section{Refinement Strategies}
\label{sec:strategies}

In this section, we discuss three strategies for the adaptive refinements in the network simulation. Errors $\Ve_m, \Ve_x, \Ve_t \in \mathbb{R}^{N_p}$ and refinements $\Vr_m, \Vr_x, \Vr_t \in \mathbb{N}^{N_p}$ are assigned to every pipe~$j \in \mathcal{J}_p$. 
One refinement in the spatial mesh or temporal discretisation is defined to be taking the new
stepsize as half of its previous value.
Hence, 
approximate relations between the initial errors $e$ and the errors after $r$ refinements $e(r)$ are given by
\begin{align}
e_{x,j} \lr{r_{x,j}} \approx  \frac{e_{x,j}}{2^{s_x \cdot r_{x,j}}}, \tab e_{t,j} \lr{r_{t,j}} \approx  \frac{e_{t,j}}{2^{s_t \cdot r_{t,j}}}, \tab \mbox{ for all\ } j \in \mathcal{J}_p, \label{discr_error_reductions}
\end{align}
where $s_x$ and $s_t$ are the convergence orders of the spatial and temporal discretisation schemes. 
The reduction
for the error is only an approximation. Therefore, to have a safe upper bound for the estimated error after refinement, we multiply these approximated errors by a \emph{safety factor of refinement} $f_r > 1$ in each of the pipes that require refinements to be made. This ensures that it is very unlikely that the error overshoots its estimated value after the refinement.

In the model hierarchy presented in Fig.~\ref{ModelH}, discretisation errors feature only in the models~\ref{M1} and~\ref{M2}. The algebraic model~\ref{M3} has no discretisation errors. Thus, when models are switched from~\ref{M3} to~\ref{M2}, discretisation errors are introduced. For pipe $j$ simulated with the most detailed model~\ref{M1}, we set the model error to $e_{m,j} =0$.

\subsection*{Strategy 1 - Individual Bounds $\lr{S1}$}

In order to meet the tolerance of the network error, we derive fixed individual error bounds for individual pipes for each of the three error types.
The simulation is then carried out such that for each pipe the error is below the individual bound for all three errors.

We set the tolerance for the model errors as
\begin{align*}
\text{tol}_m &= \kappa \cdot \text{tol}, \tab \kappa \in \lr{0,1}.
\end{align*}
The remaining tolerance is equally divided between the spatial and temporal discretisations, i.e.,
\begin{align*}
\text{tol}_x = \text{tol}_t &= \lr{1 - \kappa}/2 \cdot \text{tol}
\end{align*}
are the bounds for both error types for the entire network. To get the bounds for individual pipes, we uniformly distribute these bounds over the entire network, i.e., we divide them by the number of pipes ${N_p}$.

For the refinements, first a discretisation refinement is computed to bring the errors below the respective tolerance for each pipe. Subsequently, if the network error still exceeds the tolerance, models are refined to the next model higher up in the hierarchy. The simulation errors are then re-evaluated and this cycle is repeated until
\eqref{TolTarget} is satisfied.
This strategy is discussed in detail in~\cite{DomschkePhD}. The pseudocode for \emph{Individual Bounds} is given in Algorithm \ref{ind_bound}.

\begin{algorithm}[!b]
\caption{ individualBound}\label{ind_bound}
\begin{algorithmic}[1]
\Require $\Ve_m,\Ve_x,\Ve_t,\text{tol},s_x,s_t$
\Ensure $\Vr_m, \Vr_x, \Vr_t$
\While{networkError $ > \text{tol} \cdot M(\mathbf{u}^h)$}
	\For{$j  =  1, \ldots,  {N_p}$}
		\If{$e_{x,j}> \text{tol}_{x}\cdot M(\mathbf{u}^h)/{N_p}$ }
			\State $ r_{x,j} \gets \text{ceil} \lr{ \ds\frac{\log\lr{f_r e_{x,j} {N_p}/(\text{tol}_{x}\cdot M(\mathbf{u}^h))}}{{\log\lr{2^{s_x}}} }} $
		\EndIf
		\If{$e_{t,j}> \text{tol}_{t}\cdot M(\mathbf{u}^h)/{N_p}$ }
			\State $ r_{t,j} \gets \text{ceil} \lr{ \ds\frac{\log\lr{f_r e_{t,j} {N_p}/(\text{tol}_{t}\cdot M(\mathbf{u}^h))}}{{\log\lr{2^{s_t}}} }} $
		\EndIf
	\EndFor
	\State \textbf{update} networkError
	\If{networkError $ > \text{tol}\cdot M(\mathbf{u}^h)$}
		\For{$j = 1, \dots, {N_p}$}
			\If{$e_{m,j} > \text{tol}_{m}\cdot M(\mathbf{u}^h)/{N_p}$}
				\State $r_{m,j} \gets r_{m,j} +1 $
			\EndIf
		\EndFor
	\EndIf
	\State \textbf{update} networkError
\EndWhile
\end{algorithmic}
\end{algorithm}

\subsection*{Strategy 2 - Maximal Error Refinement $\lr{S2}$}

Since Strategy $1$ assigns individual tolerances, it loses the view of the network as a whole. However, the contribution of different errors to the overall network can balance each other without overshooting the total network tolerance. This is accounted for in the following strategy where we seek to make only those refinements which result in the maximal error reduction. This results in an iterative procedure. 
For every pipe $j$, we compute in every iteration the error reduction due to a single refinement in the model $\Delta e_{m,j}$ and in the space and time discretisations $\Delta e_{x,j}$ and $\Delta e_{t,j}$. The best option
\begin{align*}
b_j = \max_{j \in \mathcal{J}_p}\Lr{\Delta e_{m,j}, \Delta e_{x,j}, \Delta e_{t,j}}
\end{align*}
is passed to the network.  On the network level, we mark
those refinements for which the error reductions are larger than $\phi \cdot \max\Lr{b}$, with $\phi \leq 1$. This iteration is repeated until the network error is brought
below the tolerance. The function \emph{maximalErrorRefinement}, see Algorithm~\ref{max_refin}, represents the network controller and the function \emph{errorReduction}, see Algorithm~\ref{err_red}, represents the pipe level computations.

\begin{algorithm}[!b]
\caption{maximalErrorRefinement}
\label{max_refin}
\begin{algorithmic}[1]
\Require $\Vm, \Ve_m, \Ve_x, \Ve_t ,\text{tol}, s_x, s_t, \phi$
\Ensure $\Vr_m, \Vr_x, \Vr_t$
\For{$j  =  1, \ldots,  {N_p}$}
	\State $b_j, z_j\gets$ errorReduction($m_j, r_{m,j}, r_{x,j}, r_{t,j}, e_{m,j}, e_{x,j}, e_{t,j}, s_x, s_t$) \label{errorReduction_1}
\EndFor
\While {networkError $> \text{tol} \cdot M(\mathbf{u}^h)$}
	\State bound $\gets \phi\cdot\max\Lr{\Vb}$
	\For{$j  =  1, \ldots,  {N_p}$}
		\If{$b_j > $ bound}
			\State $r_{z_{j},j} \gets r_{z_{j},j} + 1 $	
			\State $b_j,z_j\gets$ errorReduction($m_j, r_{m,j}, r_{x,j}, r_{t,j}, e_{m,j}, e_{x,j}, e_{t,j}, s_x, s_t$) \label{errorReduction_2}
		\EndIf
	\EndFor
	\State \textbf{update} networkError
\EndWhile
\end{algorithmic}
\end{algorithm}

\begin{algorithm}[!b]
\caption{errorReduction}
\label{err_red}
\begin{algorithmic}[1]
\Require $m, r_m, r_x, r_t, e_m, e_x, e_t, s_x, s_t$
\Ensure $b, z$
\State \textbf{function} safety(x) \Return $ 1+ (f_r-1)\sign(\text{x})$
\State \textbf{function} spaceError(m,e,r) \Return $ \!\frac{\text{e}}{2^{s_x \cdot \text{r}}}\!\cdot\text{safety}(\text{r})\!\cdot\! F_{x}(\text{m}) $ \label{spaceError}
\State \textbf{function} timeError(m,e,r) \Return $ \!\frac{\text{e}}{2^{s_t \cdot \text{r}}}\!\cdot\text{safety}(\text{r})\!\cdot\! F_{t}(\text{m}) $ \label{timeError}
\State $m_c \gets m - r_m$
\If{$m_c \neq 1 $}
\State $\Delta e_m \gets F_m(m_c,m_c - 1) \!\cdot\! e_m$ \label{model_error_reduction} 
\Statex \hspace{1.8cm}$+\,\text{spaceError}(m_c,e_x,r_x)-\text{spaceError}(m_c-1,e_x,r_x)$
\Statex \hspace{1.8cm}$+\,\text{timeError}(m_c,e_t,r_t)-\text{timeError}(m_c-1,e_t,r_t)$
\Else
	\State $\Delta e_m \gets 0$
\EndIf
\State $\Delta e_x \gets \text{spaceError}(m_c,e_x,r_x)-\text{spaceError}(m_c,e_x,r_x+1) $
\State $\Delta e_t \gets \text{timeError}(m_c,e_t,r_t)-\text{timeError}(m_c,e_t,r_t+1) $ \label{end_main_part_alg_err_red}
\State $\LR{b ,z} = \texttt{max}\Lr{\Delta e_m, \Delta e_x, \Delta e_t}$
\end{algorithmic}
\end{algorithm}

The spatial and temporal discretisation errors also depend on the simulation model. In lines \ref{spaceError} and \ref{timeError} of Algorithm \ref{err_red}, we use $F_{x}\!\lr{m_c}$ and $F_{t}\!\lr{m_c}$, with~$m_c$ the current model, as error amplification factors for the discretisation errors which account for this model dependency. Here, $m_c\!=\!1, 2, 3$ refers to models M1, M2 and M3, respectively. Since the discretisation error is absent for $m_c\!=\!3$, we set $F_{x}\!\lr{3}\!=\!F_{t}\!\lr{3}\!=\!0$. 
Furthermore, we set $F_{x}\!\lr{2}\!=\!F_{t}\!\lr{2}\!=\!1$. For model~\ref{M1}, the amplification factor for the discretisation errors is set with respect to the benchmark model~\ref{M2}.
The factor $F_m\!\lr{a,b}$ denotes an error reduction for the model error when models are shifted from~$a$ to~$b$. The notation $\LR{b,z} = \texttt{max}\{\cdot\}$ is similar to MATLAB notation where $b$ denotes the maximal element and $z$ is the corresponding index. Note that in determining $\Delta e_m$ we also consider changes in the spatial and temporal discretisation errors. The central idea is to account for net error reduction.

\subsection*{Strategy 3 - Maximal Error-to-Cost Refinement $\lr{S3}$}

The adaptive refinements are made with an objective of reducing the computational cost without compromising on the simulation error. However, the previous two strategies do not address the computational costs explicitly. They address the error tolerance which is merely a constraint to the adaptive strategies, viewed in the optimization setting \eqref{eq:gen_Knapsack}. In Strategy 3, however,
we also take into account the computational costs that are incurred by the refinement. The idea of this strategy is 
similar to the greedy approximation algorithm for solving the unbounded knapsack problem, see \cite{Dan57}. Strategy 3, given in Algorithm~\ref{max_e2c_refin}, is similar to Algorithm~\ref{max_refin} on the network level. On the pipe level, however, we compute the cost additions $\Delta c_{m,j}, \Delta c_{x,j}, \Delta c_{t,j}$, $\mbox{ for all\ } j \in \mathcal{J}_p$, using the cost functional $F_{c}\lr{m_c, r_x, r_t}$, for each of the error reductions $\Delta e_{m,j}, \Delta e_{x,j}, \Delta e_{t,j}$. The error controller on pipe level then passes the best option
\begin{align*}
b_j = \max\Lr{\frac{\Delta e_{m,j}}{\Delta c_{m,j}},\frac{\Delta e_{x,j}}{\Delta c_{x,j}},\frac{\Delta e_{t,j}}{\Delta c_{t,j}}},
\end{align*}
i.e., the maximal error-to-cost ratio, to the network. This extension in the error controller on pipe level is given in function \emph{errorToCostReduction} in Algorithm~\ref{e2c_red}.

\begin{algorithm}[!hbt]
\caption{maximalErrorToCostRefinement}
\label{max_e2c_refin}
\begin{algorithmic}[1]
\Require $\Vm, \Ve_m, \Ve_x, \Ve_t, \text{tol}, s_x, s_t, \phi$
\Ensure $\Vr_m, \Vr_x, \Vr_t$
\setcounter{ALG@line}{1}
\Statex Same as Algorithm \ref{max_refin}, replacing lines \ref{errorReduction_1} and \ref{errorReduction_2} with
\Statex $b_j,z_j\gets$ errorToCostReduction($m, r_{m,j}, r_{x,j}, r_{t,j}, e_{m,j}, e_{x,j}, e_{t,j}, s_x, s_t$)
\end{algorithmic}
\end{algorithm}

\begin{algorithm}[!hbt]
\caption{errorToCostReduction}
\label{e2c_red}
\begin{algorithmic}[1]
\Require $m, r_m, r_x, r_t, e_m, e_x, e_t, s_x, s_t$
\Ensure $b,z$
\setcounter{ALG@line}{11}
\Statex Lines 1-\ref{end_main_part_alg_err_red} same as Algorithm \ref{err_red}
\If{$m_c \neq 1 $}
	\State $\Delta c_m \gets F_c(m_c-1,r_x,r_t) - F_c(m_c,r_x,r_t) $
\EndIf
\State $\Delta c_x \gets F_c(m_c,r_x+1,r_t) - F_c(m_c,r_x,r_t) $
\State $\Delta c_t \gets F_c(m_c,r_x,r_t+1) - F_c(m_c,r_x,r_t)$
\State $\LR{b,z} = \texttt{max}\Lr{\Delta e_m / \Delta c_m, \Delta e_x / \Delta c_x, \Delta e_t / \Delta c_t}$
\end{algorithmic}
\end{algorithm}

\section{Design of Experiment}
\label{sec:des_of_expt}

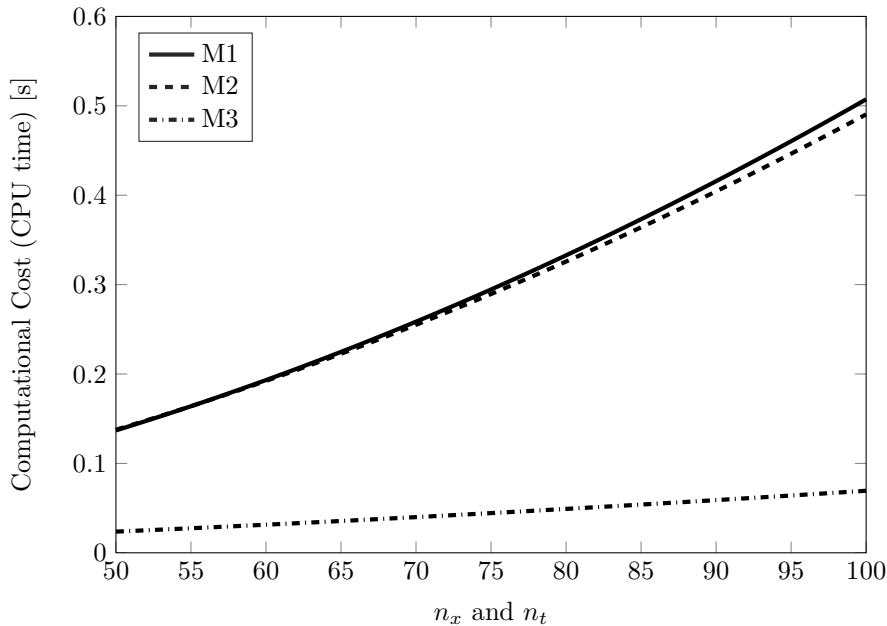
\begin{figure}[!h]
\centering
%
%
\begin{tikzpicture}

\begin{axis}[%
width=10.333in,
height=7.394in,
at={(1.733in,0.998in)},
scale only axis,
scale = 0.38,
xmin=50,
xmax=100,
xlabel={$n_x$ and $n_t$},
ymin=0,
ymax=0.6,
ylabel={Computational Cost (CPU time) [s]},
axis background/.style={fill=white},
legend style={at={(0.03,0.97)},anchor=north west,legend cell align=left,align=left,draw=white!15!black}
]
\addplot [solid,line width=1.6pt]
  table[row sep=crcr]{%
50	0.136922231422196\\
51	0.142141951190638\\
52	0.147453490020357\\
53	0.152856650736975\\
54	0.158351240333953\\
55	0.163937069808288\\
56	0.16961395400558\\
57	0.175381711473777\\
58	0.181240164324966\\
59	0.187189138104636\\
60	0.193228461667893\\
61	0.199357967062134\\
62	0.205577489415756\\
63	0.211886866832493\\
64	0.218285940290996\\
65	0.224774553549347\\
66	0.231352553054158\\
67	0.238019787854001\\
68	0.244776109516868\\
69	0.251621372051447\\
70	0.258555431831964\\
71	0.265578147526385\\
72	0.272689380027791\\
73	0.27988899238873\\
74	0.287176849758391\\
75	0.294552819322427\\
76	0.302016770245296\\
77	0.309568573614974\\
78	0.317208102389915\\
79	0.32493523134814\\
80	0.332749837038337\\
81	0.340651797732884\\
82	0.348640993382672\\
83	0.356717305573667\\
84	0.364880617485095\\
85	0.373130813849194\\
86	0.381467780912441\\
87	0.389891406398191\\
88	0.398401579470661\\
89	0.406998190700189\\
90	0.415681132029717\\
91	0.42445029674244\\
92	0.433305579430559\\
93	0.442246875965106\\
94	0.451274083466773\\
95	0.46038710027772\\
96	0.469585825934304\\
97	0.478870161140702\\
98	0.488240007743375\\
99	0.497695268706356\\
100	0.507235848087313\\
};
\addlegendentry{M1};

\addplot [dashed,line width=1.6pt]
  table[row sep=crcr]{%
50	0.137671573631479\\
51	0.142759760719257\\
52	0.147931710015856\\
53	0.153187149829151\\
54	0.158525814499418\\
55	0.163947444155273\\
56	0.169451784483784\\
57	0.175038586513684\\
58	0.180707606410719\\
59	0.186458605284266\\
60	0.192291349004397\\
61	0.198205608028679\\
62	0.204201157238026\\
63	0.210277775781003\\
64	0.216435246926\\
65	0.222673357920783\\
66	0.228991899858926\\
67	0.235390667552703\\
68	0.241869459412018\\
69	0.24842807732902\\
70	0.255066326568041\\
71	0.261784015660549\\
72	0.268580956304811\\
73	0.275456963270003\\
74	0.2824118543045\\
75	0.289445450048113\\
76	0.29655757394806\\
77	0.303748052178455\\
78	0.311016713563124\\
79	0.318363389501581\\
80	0.325787913897977\\
81	0.33329012309289\\
82	0.340869855797784\\
83	0.348526953032019\\
84	0.356261258062281\\
85	0.364072616344297\\
86	0.371960875466745\\
87	0.379925885097228\\
88	0.387967496930235\\
89	0.396085564636979\\
90	0.404279943817029\\
91	0.412550491951651\\
92	0.420897068358785\\
93	0.429319534149577\\
94	0.437817752186396\\
95	0.446391587042274\\
96	0.455040904961706\\
97	0.463765573822748\\
98	0.472565463100354\\
99	0.481440443830911\\
100	0.490390388577905\\
};
\addlegendentry{M2};

\addplot [dashdotted,line width=1.6pt]
  table[row sep=crcr]{%
50	0.0236880142404799\\
51	0.0244268047984774\\
52	0.0251736193752132\\
53	0.0259283882938387\\
54	0.0266910437928914\\
55	0.0274615199389582\\
56	0.0282397525448428\\
57	0.0290256790928012\\
58	0.0298192386624457\\
59	0.0306203718629581\\
60	0.0314290207692852\\
61	0.0322451288620204\\
62	0.0330686409707003\\
63	0.0338995032202713\\
64	0.0347376629804999\\
65	0.0355830688181228\\
66	0.0364356704515461\\
67	0.037295418707924\\
68	0.0381622654824551\\
69	0.039036163699753\\
70	0.0399170672771557\\
71	0.0408049310898509\\
72	0.0416997109377012\\
73	0.0426013635136663\\
74	0.043509846373722\\
75	0.0444251179081877\\
76	0.0453471373143765\\
77	0.0462758645704926\\
78	0.0472112604107008\\
79	0.0481532863013037\\
80	0.0491019044179623\\
81	0.0500570776239022\\
82	0.0510187694490511\\
83	0.0519869440700576\\
84	0.0529615662911416\\
85	0.0539426015257358\\
86	0.054930015778873\\
87	0.0559237756302826\\
88	0.05692384821816\\
89	0.0579302012235738\\
90	0.0589428028554785\\
91	0.0599616218363052\\
92	0.0609866273880982\\
93	0.0620177892191745\\
94	0.0630550775112784\\
95	0.0640984629072098\\
96	0.0651479164989024\\
97	0.0662034098159319\\
98	0.0672649148144338\\
99	0.0683324038664126\\
100	0.069405849749424\\
};
\addlegendentry{M3};

\end{axis}
\end{tikzpicture}%
\caption{Computational costs for the simulation of a single pipe using different models from the hierarchy against the number of nodes in space $n_x$ and time $n_t$.}
\label{plot:CostFunc}
\end{figure}

We tested the performance of the three refinement strategies on $10^4$ random samples of error distributions. The random samples represent errors in a simulation of gas flow in a given network of ${N_p} = 12$ pipes. A better strategy will lead to lower
computational costs while reducing the simulation error to a level below the tolerance. We compute the computational cost per pipe using a cost functional of the form
\begin{align}\label{costFunctionals}
F\lr{m,n_x,n_t} = C_m \cdot n_x^{\alpha_m} \cdot n_t^{\beta_m},
\end{align}
where $m \in \{1,2,3\}$ denotes the model and $C_m$, $\alpha_m$ and $\beta_m$ are model-dependent constants. 
These constants, given in Table~\ref{CostFunc}, are determined by the method of least squares. 
For this, gas flow simulations through a single pipe are performed using the software \emph{ANACONDA} (cf.\ \cite{Kolb2011,Kolb2009}) with many different values of~$n_x$ and~$n_t$, which return the corresponding computational cost values~$F$. A plot showing the computational costs of simulating a single pipe using different models depending on the number of nodes in space and time is depicted in Fig.~\ref{plot:CostFunc}. We can rewrite the functional~(\ref{costFunctionals}) in terms of refinements assuming that the initial number of nodes $n_{x,0}, n_{t,0}$ are known. Then we get
\begin{align}\label{costFunctionals_2}
F_c\lr{m,r_x,r_t} = C_m  \lr{2^{r_x} \cdot n_{x,0}}^{\alpha_m} \cdot \lr{2^{r_t} \cdot n_{t,0}}^{\beta_m}.
\end{align}

\begin{table}[!b]
\centering
\caption{Cost functional constants in \eqref{costFunctionals} and \eqref{costFunctionals_2}.}
\label{CostFunc}
\vspace{0.2cm}
\begin{tabular}{cccc}
\hline
\textbf{$m$}    & \textbf{$C_m$}  & \textbf{$\alpha_m$} & \textbf{$\beta_m$} \\ \hline
$1$ & \num{8.45e-5} & $0.952$          & $0.937$         \\ 
$2$ & \num{1.06e-4} & $0.908$          & $0.925$         \\ 
$3$ & \num{5.49e-5} & $0.694$          & $0.857$         \\ \hline
\end{tabular}
\end{table}

The initial number of space and time discretisation nodes are chosen from the interval $\LR{100, 200}$ using a random number generator. All models are set to the most simple model~\ref{M3} in the beginning.
The error reduction upon model refinement also takes into account the introduction or increase of the spatial and temporal errors. This requires that the spatial and temporal errors are small when compared to the model error. Hence, for the experiment, the initial model errors are taken from $\mathcal{U}\LR{ 0, 1 }$, where $\mathcal{U}\LR{ a, b }$ denotes a uniform probability distribution on the interval $\LR{a,b}$, and the spatial and temporal discretisation errors are taken from $\mathcal{U}\LR{ 0, 0.2 }$. 
For the model and discretisation errors after refinement, the approximate error reductions in lines \ref{spaceError}, \ref{timeError} and \ref{model_error_reduction} of Algorithm~\ref{err_red} are used, where we choose $F_m(3,2) \! = \! 3/4$ and $F_m(2,1) \! = \! 1/4$.
Models~\ref{M1} and~\ref{M2} are discretised with the implicit box scheme, as in \cite{KolLB10}, for which it holds that $s_x \! = \! 2$ and $s_t \! = \! 1$.
The parameter $\kappa = 1/3$ for Strategy~1 is chosen such that all three errors have an equal fraction of the tolerance. Strategies~2 and~3 are tested for a fraction $\phi \in \Lr{0.8, 0.9, 1}$ of the maximal best option. The strategies work for a relative error tolerance of $\text{tol} = \num{e-1}$ with a target functional value $M(\mathbf{u}^h) = 2.5 \cdot {N_p} = 30$. Hence, we require for the total network simulation error that
\begin{align*}
\sum_{j \in \mathcal{J}_p} \lr{ e_{m,j} + e_{x,j} + e_{t,j} }< \text{tol} \cdot M(\mathbf{u}^h) = 3.
\end{align*}
The results of this experiment are given in the next section.

\section{Results and Discussions}
\label{sec:results}

Each strategy returns a refinement scheme which brings the simulation error below the tolerance. The goal is to have low computational costs. 
The mean of the total computational cost values in CPU seconds over $10^4$ samples is shown in Table~\ref{samp10k}. 
We show the percentage savings in mean total computational cost of the strategies with respect to Strategy $1$.
We denote the strategies as $S1 - S3$. The subscripts $1,2,3$ refer to $\phi = 0.8,0.9$ and $1$, respectively. We observe that strategies $S2$ and $S3$ have a percentage saving of over $77\%$
with respect to $S1$ for all values of $\phi \in \Lr{0.8, 0.9, 1}$. Among the different values, $\phi = 1$ performs best for both $S2$ and $S3$. 

Thus, by working with a greedy-like strategy for error control, an adaptive process can reduce the computational cost significantly. Furthermore, accounting for the computational cost explicitly in our estimates, we find even better refinement schemes that result in lower computational costs.

\begin{table}[!h]
\centering
\caption{Mean total cost values in CPU seconds for strategies $S1$ - $S3$ and savings with respect to $S1$. The subscripts $1,2,3$ denote the different values of $\phi = 0.8, 0.9,1$, respectively. }
\label{samp10k}
\vspace{0.3cm}
\begin{tabular}{cccccccc}
\hline
\textbf{strategies}              & $S1$ & $S2_1$ & $S2_2$ & $S2_3$ & $S3_1$ & $S3_2$ & $S3_3$ \\ \hline
\textbf{Mean Total Cost}      & $36.1$    & $8.30$      & $8.04$      & $7.72$      & $7.64$      & $7.39$      & $7.12$      \\ 
\textbf{Savings (\%)} & -         & $77.0$      & $77.7$      & $78.6$      & $78.8$      & $79.5$      & $80.3$      \\ \hline
\end{tabular}
\end{table}

\section{Application to a Realistic Network Simulation}
\label{Application_to_a_Network_Simulation}

\begin{figure}[!bp]
\centering \small
 \begin{tikzpicture}[scale=0.088]
		\node[circle,fill=black,draw,scale=0.3] (N01) at (20,40) {};
		\node[circle,fill=black,draw,scale=0.3] (N04) at (30,40) {};
		\node[circle,fill=black,draw,scale=0.3] (N05) at (40,40) {};
		\node[circle,fill=black,draw,scale=0.3] (N06) at (50,40) {};
		\node[circle,fill=black,draw,scale=0.3] (N08) at (60,40) {};
		\node[circle,fill=black,draw,scale=0.3] (N09) at (70,40) {};
		\node[circle,fill=black,draw,scale=0.3] (N10) at (90,60) {};
		\node[circle,fill=black,draw,scale=0.3] (N11) at (90,20) {};
		\node[circle,fill=black,draw,scale=0.3] (N12) at (90,0) {};
		\node[circle,fill=black,draw,scale=0.3] (N13) at (110,40) {};
		\node[circle,fill=black,draw,scale=0.3] (N16) at (120,40) {};
    	   	\node[circle,draw,scale=0.8,label={[scale=1]left:{S01}}] (S01) at (10,40) {};
    	   	\node[circle,draw,scale=0.8,label={[scale=1]above:{S02}}] (S02) at (20,45) {};
 		\node[draw,scale=1,label={[scale=1]above:{C01}}] (C01) at (50,60) {};
 		\node[draw,scale=1,label={[scale=1]above:{C02}}] (C02) at (110,0) {};
 		\node[draw,scale=1,label={[scale=1]above:{C03}}] (C03) at (130,50) {};
 		\node[draw,scale=1,label={[scale=1]below:{C04}}] (C04) at (130,30) {};
 		\path[-latex]
 		(S01) edge node[above,sloped,scale=1,text=black] {{P01}} (N01)
 		(N01) edge node[above,sloped,scale=1,text=black] {{P02}} (N04)
 		(N05) edge node[above,sloped,scale=1,text=black] {{P03}} (N06)
 		(N06) edge node[above,sloped,scale=1,text=black] {{P04}} (C01)
 		(N06) edge node[above,sloped,scale=1,text=black] {{P05}} (N08)
 		(N09) edge node[above,sloped,scale=1,text=black] {{P06}} (N10)
 		(N10) edge node[above,sloped,scale=1,text=black] {{P07}} (N13)
 		(N09) edge node[above,sloped,scale=1,text=black] {{P08}} (N11)
 		(N11) edge node[above,sloped,scale=1,text=black] {{P09}} (N13)
 		(N12) edge node[above,sloped,scale=1,text=black] {{P10}} (C02)
 		(N16) edge node[below,sloped,scale=1,text=black] {{P11}} (C03)
 		(N16) edge node[below,sloped,scale=1,text=black] {{P12}} (C04)
 		(S02) edge (N01)
		(N11) edge node[above,sloped] {{CV01}} (N12)
		(N11) edge node[left=4] {} (N12)
 		;
		\draw[fill=black] (N11)++(0,-10)++(0,-1) -- ++(0.5,2) -- ++(-1,0) -- cycle;
		\path[draw] (N04) -- ++(2,2.5) -- node[above=0.1] {{Comp01}} ++(6,0) -- ++(2,-2.5) -- ++(-2,-2.5) --
		node[below=4] {} ++(-6,0) -- ++(-2,2.5);
		\draw[fill=white] (N04)+(5,2.5) circle (2);
		\path[draw] (N04)++(5,4.5) -- ++(2,-2) -- ++(-2,-2);
		\draw[fill=black] (N04)++(5,-2.5)++(-1,0.5) -- ++(2,-1) -- ++(0,1) -- ++(-2,-1) -- cycle;
		\path[draw] (N08) -- ++(2,2.5) -- node[above=0.1] {{Comp02}} ++(6,0) -- ++(2,-2.5) -- ++(-2,-2.5) --
		node[below=4] {} ++(-6,0) -- ++(-2,2.5);
		\draw[fill=white] (N08)+(5,2.5) circle (2);
		\path[draw] (N08)++(5,4.5) -- ++(2,-2) -- ++(-2,-2);
		\draw[fill=black] (N08)++(5,-2.5)++(-1,0.5) -- ++(2,-1) -- ++(0,1) -- ++(-2,-1) -- cycle;
		\path[draw] (N13) -- ++(2,2.5) -- node[above=0.1] {{Comp03}} ++(6,0) -- ++(2,-2.5) -- ++(-2,-2.5) --
		node[below=4] {} ++(-6,0) -- ++(-2,2.5);
		\draw[fill=white] (N13)+(5,2.5) circle (2);
		\path[draw] (N13)++(5,4.5) -- ++(2,-2) -- ++(-2,-2);
		\draw[fill=black] (N13)++(5,-2.5)++(-1,0.5) -- ++(2,-1) -- ++(0,1) -- ++(-2,-1) -- cycle;
 \end{tikzpicture}
\caption{Gas supply network with compressor stations and a control valve.}
\label{fig:TestNet3}
\end{figure}
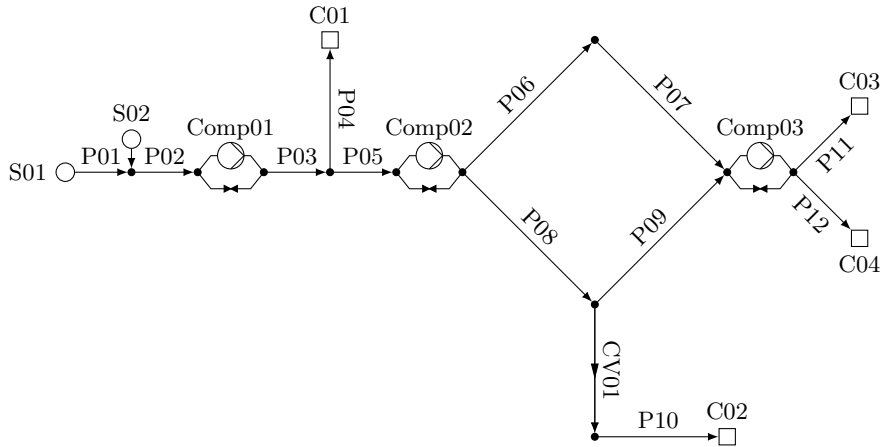

We now apply the three different strategies to a simulation of a gas supply network, which is shown in Fig.~\ref{fig:TestNet3}. The considered network consists of twelve pipes (P01 -- P12, with lengths between 30km and 100km), two sources (S01 -- S02), four consumers (C01 -- C04), three compressor stations (Comp01 - Comp03) and one control valve (CV01). Starting with stationary initial data, the boundary conditions and the control for the compressor stations and the control valve are time-dependent. The simulation time is \SI{14400}{\second}. The target functional $M(\mathbf{u})$ is given by the total fuel gas consumption of the three compressors and the error estimators are evaluated using a dual weighted residual method as developed, for example, in  \cite{Becker2001,Braack2003,Rannacher2009,Besier2012}. 
For details of the derivation of the error estimators we refer to \cite{DomschkePhD,Domschke2015}.
The simulation is performed using the software \emph{ANACONDA} (cf.\ \cite{Kolb2011,Kolb2009}).

\begin{remark}{\rm 
For the
strategies proposed in Section \ref{sec:strategies}, the temporal error was considered individually for each pipe. Since in the implementation of \emph{ANACONDA}, the time stepping is uniform for the entire network, the temporal error~$\eta_t$ was computed globally and divided by the number of pipes in order to get a local temporal error. However, if a best option $B_j$ was supposed to be the temporal error and had to be refined then all pipes were refined in time uniformly and the best options $B_j$ were updated.}
\end{remark}

A reference solution was computed using model \ref{M1} and a very fine discretisation. The strategies were run with a relative tolerance of $\text{tol}=\num{e-4}$.
Table~\ref{tab:resultsExample} shows the relative error of the simulation compared to the reference solution, the CPU time taken and the percentage savings of strategy $S2$ and $S3$ in relation to strategy $S1$.

Compared to the synthetic experiment in Section~\ref{sec:results}, we see that the savings of strategies~$S2$ and~$S3$ applied to the simulation are in a similar range. The choice of the parameter $\phi\in\{0.8,0.9,1.0\}$, however, does not seem to have a significant influence on the saving. Moreover, the \emph{Maximal Error-to-Cost} strategy $S3$ does not result in a larger saving of CPU time. What is noticeable is that the relative errors of the strategies $S2$ and $S3$ are closer to the proposed relative tolerance, which shows that they are not as restrictive as the \emph{Individual Bounds} strategy~$S1$.

\begin{table}[!ht]
\centering
\caption{Relative error ($\text{tol}=\num{e-4}$) and computational costs (in \si{s}) of a network simulation using strategies $S1$ - $S3$ and savings with respect to $S1$. The subscripts $1,2,3$ denote the different values of $\phi \in \{0.8, 0.9,1\}$, respectively. }
\label{tab:resultsExample}
\vspace{0.3cm}
\begin{tabular}{cccc}
 \hline
 Strategy & Relative Error & CPU time [s] & Savings\\
 \hline
 $S1_{~}$ & \num{1.384e-5} & 7.53 & -\\
 $S2_1$ & \num{3.730e-5} & 2.42 & \SI{67.9}{\percent}\\
 $S2_2$ & \num{3.091e-5} & 2.27 & \SI{69.9}{\percent}\\
 $S2_3$ & \num{3.091e-5} & 2.29 & \SI{69.9}{\percent}\\
 $S3_1$ & \num{5.076e-5} & 3.12 & \SI{58.6}{\percent}\\
 $S3_2$ & \num{5.190e-5} & 2.83 & \SI{62.4}{\percent}\\
 $S3_3$ & \num{5.098e-5} & 2.83 & \SI{62.4}{\percent}\\
 reference   & - & 1013 \\
 \hline
\end{tabular}
\end{table}

\section{Conclusions}
\label{Conclusions}

In this paper we address the problem of automatic error control for large scale gas flow simulations that use a model hierarchy.
The simulation needs to be reliable, i.e., keeping the total relative error below a specified tolerance, while retaining low computational costs.
The problem of finding an optimal refinement strategy is a generalisation of the knapsack problem.
We present three strategies for adaptive simulation error control via spatial and temporal discretisation mesh and model refinements.
The strategy \emph{Individual Bounds}, which is currently implemented in \emph{ANACONDA}, sets a uniform bound for each error type and each pipe, \emph{Maximal Error Refinement} iteratively chooses those refinements that result in the largest error reduction and has a network overview, and \emph{Maximal Error-to-Cost Refinement} also accounts for the increase in computational cost inflicted by the refinement.

We constructed
a synthetic experiment to test the three strategies. From this experiment we
see that the two greedy-like strategies significantly reduce the computational cost as compared to the \emph{Individual Bounds} strategy.
This result is largely reflected in an actual gas flow simulation using \emph{ANACONDA} for a~$12$ pipe network including compressor stations and a control valve. Especially when the simulation process is a key component in a gas flow optimisation problem, the novel refinement strategies lead to considerable computational savings without compromising on the simulation accuracy.

\section*{Acknowledgements}
The authors would like to thank the Deutsche Forschungsgemeinschaft for their support within the projects B01 and B03 of the Collaborative Research Centre\slash{}Trans\-regio~154 \textit{Mathematical Modelling, Simulation and Optimization using the Example of Gas Networks}. Jens Lang was also supported by the Excellence Initiative of the German Federal and State Governments via the Darmstadt Graduate School of Excellence Energy Science and Engineering.

\section*{References}

\bibliography{AMC-literature}

\end{document}